\documentclass[a4paper,12pt]{article}
\usepackage{amsmath, amssymb,amsfonts}
\usepackage[all]{xy}
\makeatletter
\newbox\bk@bxb
\newbox\bk@bxa
\newif\if@bkcont
\newcount\bk@lcnt

\def\breakboxskip{2pt}
\def\breakboxparindent{1.8em}

\def\breakbox{\vskip\breakboxskip\relax
\setbox\bk@bxb\vbox\bgroup
\advance\linewidth -2\fboxrule
\hsize\linewidth\@parboxrestore
\parindent\breakboxparindent\relax}

\def\bk@split{%
\@tempdimb\ht\bk@bxb 
\advance\@tempdimb\dp\bk@bxb
\setbox\bk@bxa\vsplit\bk@bxb to\z@ 
\setbox\bk@bxa\vbox{\unvbox\bk@bxa}
\setbox\@tempboxa\vbox{\copy\bk@bxa\copy\bk@bxb}
\advance\@tempdimb-\ht\@tempboxa
\advance\@tempdimb-\dp\@tempboxa}

\def\bk@addfsepht{%
\setbox\bk@bxa\vbox{\vskip\fboxsep\box\bk@bxa}}

\def\bk@addskipht{%
\setbox\bk@bxa\vbox{\vskip\@tempdimb\box\bk@bxa}}

\def\bk@addfsepdp{%
\@tempdima\dp\bk@bxa
\advance\@tempdima\fboxsep
\dp\bk@bxa\@tempdima}

\def\bk@addskipdp{%
\@tempdima\dp\bk@bxa
\advance\@tempdima\@tempdimb
\dp\bk@bxa\@tempdima}

\def\bk@line{%
\hbox to \linewidth{%
\hskip-2\fboxsep\vrule \@width\fboxrule\hskip.5\fboxsep\vrule \@width\fboxrule\hskip1.5\fboxsep
\box\bk@bxa\hfil
}}%

\def\endbreakbox{\egroup
\ifhmode\par\fi{\noindent\bk@lcnt\@ne
\@bkconttrue\baselineskip\z@\lineskiplimit\z@
\lineskip\z@\vfuzz\maxdimen
\bk@split\bk@addfsepht\bk@addskipdp
\ifvoid\bk@bxb 
\def\bk@fstln{\bk@addfsepdp
\hskip-\parindent\vbox{\llap{\raisebox{-2ex}{\rule{1.5\fboxsep}{\fboxrule}\hskip.5\fboxsep}}\bk@line\llap{\rule{1.5\fboxsep}{\fboxrule}\hskip.5\fboxsep}}}

\else 
\def\bk@fstln{\vbox{\llap{\raisebox{-2ex}{\rule{1.5\fboxsep}{\fboxrule}\hskip.5\fboxsep}}\bk@line}\hfil%
\advance\bk@lcnt\@ne
\loop
\bk@split\bk@addskipdp\leavevmode
\ifvoid\bk@bxb 
\@bkcontfalse\bk@addfsepdp
\vtop{\bk@line\llap{\rule{2\fboxsep}{\fboxrule}}}%

\else 
\bk@line
\fi
\hfil\advance\bk@lcnt\@ne
\if@bkcont\repeat}%
\fi
\leavevmode\bk@fstln\par}\vskip\breakboxskip\relax}

\def\smp{\smallskip\par}

\def\pf{\noindent{\bf Proof~:}\ }

\def\findemo{~\leaders\hbox to 1em{\hss\  \hss}\hfill~\raisebox{.5ex}{\framebox[1ex]{}}\smp}

\def\mpn{\medskip\par\noindent}
\def\smpn{\smallskip\par\noindent}

\def\smp{\smallskip\par}
\def\smpn{\smallskip\par\noindent}

\def\mpoint{\;\;.}
\def\mvirg{\;\;,}

\def\Id{{\rm Id}}

\def\op{^{op}}

\def\N{\mathbb{N}}

\def\R{\mathbb{R}}

\newcommand{\romain}[1]{\uppercase\expandafter{\romannumeral #1}}

\newcommand{\flh}[2]{\mathop{\hbox to 12mm{\rightarrowfill}}_{\displaystyle #2}^{\displaystyle #1}\limits}
\newcommand{\sflh}[2]{\mathop{\hbox to 12mm{\rightarrowfill}}_{\scriptstyle #2}^{\scriptstyle #1}\limits}

\newcommand{\sur}[1]{\,\overline{\! #1}}

\def\op{^{op}}

\newcommand{\carre}[8]{\begin{array}{ccc}
#1&\mathop{\hbox to 12mm{\rightarrowfill}}^{\displaystyle{#2}}\limits&#3\\
\llap{$\displaystyle{#4}$}\left\downarrow\vbox to 6mm{}\right. & & \left\downarrow\vbox to 6mm{}\right.\rlap{$\displaystyle{#5}$}\\
#6&\mathop{\hbox to 12mm{\rightarrowfill}}_{\displaystyle #7}\limits&#8\\
\end{array}}

\newcommand{\carrem}[8]{\begin{array}{ccc}
#1&\mathop{\hbox to 12mm{\rightarrowfill}}^{\displaystyle #2}\limits&#3\\
\llap{$\displaystyle #4$}\left\uparrow\vbox to 6mm{}\right. & & \left\uparrow\vbox to 6mm{}\right.\rlap{$\displaystyle #5$}\\
#6&\mathop{\hbox to 12mm{\rightarrowfill}}_{\displaystyle #7}\limits&#8\\
\end{array}}

\newenvironment{enonce}[1]{\pagebreak[2]\refstepcounter{subsection}\refstepcounter{prop}\smpn{{\bf \thesection.\arabic{prop}.\ \ #1~:}}\begin{it} }{\end{it}\smp}
\newenvironment{enonce*}[1]{\pagebreak[2]\smpn{#1~:}\begin{it} }{\end{it}\smp}
\newcommand{\result}[1]{\begin{enonce}{#1}}
\def\fresult{\end{enonce}}
\newcommand{\npar}{\smallskip\par\noindent\pagebreak[2]\refstepcounter{subsection}\refstepcounter{prop}{\bf \thesection.\arabic{prop}.\ \ }}

\newenvironment{mth}[1]{\begin{breakbox}\begin{enonce}{#1}}{\end{enonce}\end{breakbox}}
\newenvironment{mth*}[1]{\begin{breakbox}\begin{enonce*}{#1}}{\end{enonce*}\end{breakbox}}

\makeatletter
\renewenvironment{enumerate}{\ifnum \@enumdepth >3 \@toodeep\else
      \advance\@enumdepth \@ne
      \edef\@enumctr{enum\romannumeral\the\@enumdepth}\list
      {\csname label\@enumctr\endcsname}{\setlength{\topsep}{1ex}\setlength{\itemsep}{0pt}\usecounter
        {\@enumctr}\def\makelabel##1{\hss\llap{##1}}}\fi}{\endlist}
\renewenvironment{itemize}{\ifnum \@itemdepth >3 \@toodeep\else \advance\@itemdepth \@ne
\edef\@itemitem{labelitem\romannumeral\the\@itemdepth}%
\list{\csname\@itemitem\endcsname}{\setlength{\topsep}{1ex}\setlength{\itemsep}{0pt}\def\makelabel##1{\hss\llap{##1}}}\fi}
{\endlist}
\makeatother
\makeatletter
\def\@sect#1#2#3#4#5#6[#7]#8{\ifnum #2>\c@secnumdepth
    \let\@svsec\@empty\else
    \refstepcounter{#1}\edef\@svsec{\csname the#1\endcsname .\hskip .5em}\fi
    \@tempskipa #5\relax
     \ifdim \@tempskipa>\z@
       \begingroup #6\relax
         \@hangfrom{\hskip #3\relax\@svsec}{\interlinepenalty \@M #8\par}%
       \endgroup
      \csname #1mark\endcsname{#7}\addcontentsline
        {toc}{#1}{\ifnum #2>\c@secnumdepth \else
                     \protect\numberline{\csname the#1\endcsname}\fi
                   #7}\else
       \def\@svsechd{#6\hskip #3\relax  
                  \@svsec #8\csname #1mark\endcsname
                     {#7}\addcontentsline
                          {toc}{#1}{\ifnum #2>\c@secnumdepth \else
                            \protect\numberline{\csname the#1\endcsname}\fi
                      #7}}\fi
    \@xsect{#5}}
\def\section{\@startsection {section}{1}{\z@}{-3.5ex plus-1ex minus
    -.2ex}{2.3ex plus.2ex}{\reset@font\Large\bf}}  

\makeatother
\renewenvironment{equation}{\refstepcounter{subsection}\refstepcounter{prop}$$}{\leqno{\bf (\theprop)}$$}

\def\mar[#1]{\ar@{-}[#1]|-{\object@{<}}}
\def\marb[#1]{\ar@{-}[#1]|{\object+{  }}}

\newenvironment{maliste}%
{ \begin{list}%
        {-}%
        {\setlength{\labelwidth}{4ex}%
         \setlength{\leftmargin}{2ex}%
	\setlength{\parskip}{0pt}%
	\setlength{\parsep}{0pt}%
	\setlength{\itemindent}{0pt}}}%
{ \end{list} }
\def\bmid{\,\big|\,}
\begin{document}
\centerline{\LARGE\bf The poset of posets}\vspace{2ex}
\centerline{\bf Serge Bouc}\vspace{6ex}
{\footnotesize {\bf Abstract:} Let $X$ be a finite set. This paper describes some topological and combinatorial properties of the poset $\Omega_X$ of order relations on $X$. In particular, the homotopy type of all the intervals in $\Omega_X$ is precisely determined, and the M\"obius function of $\Omega_X$ is computed.
}
\vspace{1ex}\par
{\footnotesize {\bf AMS subject classification: 06A11, 05E45, 55U10} }
\vspace{1ex}\par
{\footnotesize {\bf Keywords: poset, order, M\"obius function, homotopy type}}
\section{Introduction}
Let $X$ be a finite set. In this paper, we consider some topological and combinatorial properties of the {\em poset of posets} on $X$: it is the set $\Omega_X$ of order relations on $X$, ordered by inclusion of subsets of the cartesian product $X\times X$. \par
In particular we show that the the intervals $]R,S[_{\Omega_X}$ in this poset, for $R\subseteq S$, are either contractible, or have the homotopy type of a sphere of dimension $|S-R-2|$. As a consequence, we determine the M\"obius function of the poset $\Omega_X$. We also show that the upper intervals $]R,\,.\,[_{\Omega_X}$ are either contractible, or have the homotopy type of a sphere $S^{d_R}$. Each case can be precisely determined, and in particular the dimension $d_R$ of the sphere involved in the second case can be computed explicitly from the relation $R$.\par
The initial motivation for considering the poset $\Omega_X$ is a joint work with Jacques Th\'evenaz (\cite{both-relations}), in which the M\"obius function of $\Omega_X$ appears (even though the exact value of this function was not needed for our purpose in that paper). Apart from methods introduced by Quillen in his seminal paper~\cite{quillen}, the present paper is self contained.
\section{Relations, posets}
Let $X$ denote a finite set, of cardinality $n$. 
\begin{itemize}
\item A {\em relation} on $X$ is by definition a subset of the cartesian product $X\times X$. 
\item A {\em subrelation} of a relation $S$ is just a subset $R\subseteq S$ of $S$. A {\em proper} subrelation of $S$ is a proper subset $R\subset S$.
\item If $R$ and $S$ are relations on $X$, their {\em composition} $R\circ S$, also called the {\em product} $RS$ of $R$ and $S$, is defined by
$$R\circ S=\{(x,y)\in X\times X\mid\exists z\in X,\;(x,z)\in R\;\hbox{and}\;(z,y)\in S\}\mpoint$$
\item The relation $\Delta(=\Delta_X)=\{(x,x)\mid x\in X\}$ is called the {\em equality relation} on $X$. \par
\item If $R$ is a relation on $X$, then the {\em opposite} relation $R\op$ is defined by
$$R\op=\{(x,y)\in X\times X\mid (y,x)\in R\}\mpoint$$
\item A relation $R\subseteq X\times X$ is called:
\begin{maliste}
\item {\em reflexive} if $(x,x)\in R,\;\forall x\in X$. Equivalently $\Delta\subseteq R$. 
\item {\em transitive} if $\forall x,y,z\in X, \;(x,y)\in R\; \hbox{and}\; (y,z)\in R \Rightarrow (x,z)\in R$. Equivalently $R^2\subseteq R$.
\item {\em a preorder} if $R$ is reflexive and transitive. Equivalently $\Delta\subseteq\nolinebreak R=\nolinebreak R^2$.
\item {\em symmetric} if $(x,y)\in R \Rightarrow (y,x)\in R$. Equivalently $R=R\op$.
\item {\em an equivalence relation} if it is a preorder, and $R$ is symmetric. Equivalently $\Delta\subseteq R=R\op=R^2$.
\item {\em antisymmetric} if $\forall x,y\in X, \;(x,y)\in R\;\hbox{and}\;(y,x)\in R\Rightarrow x=y$. Equivalently $R\cap R\op\subseteq \Delta$. A subrelation of an antisymmetric relation is antisymmetric.
\item {\em an order} if $R$ is a preorder, and $R$ is antisymmetric. Equivalently $R=R^2$ and $R\cap R\op=\Delta$.
\end{maliste}
\item A set endowed with an order relation is called a partially ordered set, or {\em poset} for short.
\end{itemize}
Let $\Omega_X$ denote the set of order relations on $X$. The inclusion of subsets of $X\times X$ induces an order relation on $\Omega_X$. The poset $(\Omega_X,\subseteq)$ is called  the {\em poset of posets} on $X$. 
\section{The Frattini subrelation}
From now on, the set $X$ will be fixed, and often understood. In particular, the poset $\Omega_X$ will be denoted by $\Omega$.
\par
\pagebreak[4]
\begin{mth}{Notation} \label{intervals}Let $R$ be an order on $X$. For $x,y\in X$, write $x\leq_Ry$ if $(x,y)\in R$, and $x<_Ry$ if $(x,y)\in R$ and $x\neq y$. Set moreover
\begin{eqnarray*}
{[x,y]_R}&=&\{z\in X\mid x\leq_Rz\leq_R y\}\mvirg\\
{]x,y[_R}&=&\{z\in X\mid x<_Rz<_R y\}\mvirg\\
{]\,.\,,x]_R}&=&\{z\in X\mid z\leq_Rx\}\\
{]\,.\,,x[_R}&=&\{z\in X\mid z<_Rx\}\\
{[x,\,.\,[_R}&=&\{z\in X\mid x\leq_Rz\}\\
{]x,\,.\,[_R}&=&\{z\in X\mid x<_Rz\}\mpoint\\
\end{eqnarray*}
If $R\in\Omega $, a pair $(x,y)\in X\times X$ is called {\em adjacent} for $R$ if $x<_Ry$ and $]x,y[_R=\emptyset$.
The set of adjacent pairs for the relation $R$ is denoted by $\mathcal{M}_R$.
\end{mth}
\npar Let $R$ be  reflexive relation on $X$. Then $\Delta\subseteq R\subseteq R^2\subseteq \ldots\subseteq R^m$, and there is an integer $m$ such that $R^m=R^{m+1}$. This limit $\sur{R}=R^m$ is a preorder, called the {\em transitive closure} of $R$.\par
Let $R,S\in \Omega $. Notation~\ref{intervals} becomes
\begin{eqnarray*}
{[R,S]_\Omega}&=&\{T\in\Omega \mid R\subseteq T\subseteq S\}\\
{]R,S[_\Omega}&=&\{T\in\Omega \mid R\subset T\subset S\}\\
{]R,\,.\,[_\Omega}&=&\{T\in\Omega \mid R\subset T\}\mpoint
\end{eqnarray*}
An order $R$ is said to be {\em maximal} in an order $S$ if the pair $(R,S)$ is adjacent in the poset $\Omega$, i.e. if $R\subset S$ and $]R,S[_\Omega=\emptyset$.\par
The subposet $[R,S]_\Omega$ of $\Omega $ is a {\em lattice}: the meet  $T\wedge T'$ of $T,T'\in [R,S]_\Omega$ is the intersection $T\cap T'$, and their join is $T\vee T'=\sur{T\cup T'}$ (which is indeed an order, as it is a preorder contained in $S$).
\begin{mth}{Notation} When $S$ is a relation on $X$, denote by $S^{(2)}=S\times S\op$ the relation on  $X\times X$ defined by
$$\big((x,x'),(y,y')\big)\in S^{(2)}\Leftrightarrow (x,y)\in S\;\hbox{and}\;(y',x')\in S\mpoint$$
If $S$ is an order on $X$, then $S^{(2)}$ is an order on $X\times X$. In this case
$$(x,y)\leq_{S^{(2)}}(x',y')\;\Leftrightarrow\;x\leq_Sx' \;\hbox{and}\; y'\leq_Sy\mpoint$$
\end{mth}
\begin{mth}{Lemma} \label{maximal}\begin{enumerate} 
\item Let $S$ be an order on $X$, and $(x,y)$ be an adjacent pair for $S$. Then $R=S-\{(x,y)\}$ is an order on $X$, and $R$ is maximal in $S$.
\item Let $R,S\in \Omega $ with $R\subset S$. If $R$ is maximal in $S$, then 
the set $S-R$ consists of a single pair $(x,y)$, which is adjacent for $S$, and a minimal element of $(X\times X)-(R\cup R\op)$ for the relation $R^{(2)}$.
\item Let $R\in\Omega $, and let $(x,y)$ be a minimal element of $(X\times X)-(R\cup R\op)$ for the relation $R^{(2)}$. Then $S=R\sqcup\{(x,y)\}$ is an order on $X$, and $R$ is maximal in $S$.
\end{enumerate}
\end{mth}
\pf For Assertion 1, the relation $R=S-\{(x,y)\}$ is reflexive, as $x\neq y$. It is also antisymmetric, since it is contained in $S$. Showing that it is transitive amounts to showing that if $a,b,c\in X$ with $(a,b)\in R$ and $(b,c)\in R$, then $(a,c)\in R$. The assumptions imply $(a,c)\in S$, as $S$ is transitive. If $(a,c)\notin R$, then $(a,c)=(x,y)$. Then $a=x$ and $c=y$, and moreover $(x,b)\in S$ and $(b,y)\in S$. Since $(x,y)$ is adjacent for $S$, it follows that $x=b$ or $b=y$. If $x=b$, then $(b,c)=(x,y)\notin R$, and if $b=y$, then $(a,b)=(x,y)\notin R$. this contradiction shows that $(a,c)\in R$, and Assertion 1 follows.\par
For Assertion 2, let $R$ be maximal in $S$, and let $(x,y)\in S-R$ (so that in particular $x\neq y$). Then $(x,y)\in (X\times X)-(R\cup R\op)$, for otherwise $(x,y)$ and $(y,x)$ are both in $S$. Moreover, the transitive closure $\sur{R\sqcup\{(x,y)\}}$ is equal to~$S$. Thus for any $(a,b)\in S-R$, there exists a chain of pairs 
\begin{equation}\label{sequence}
(a_0,a_1), (a_1,a_2),\ldots,(a_{k-1},a_k)
\end{equation}
with $a_0=a$ and $a_k=b$, and $(a_i,a_{i+1})\in R\sqcup\{(x,y)\}$, for $i\in\{0,\ldots,k-1\}$. Since $R$ is transitive, we can assume that there are no two consecutive $i$ in $\{0,\ldots,k-1\}$ such that $(a_i,a_{i+1})\in R$. If $(a_i,a_{i+1})=(x,y)$ for two distinct values $l$ and $l'$ of $i$, then we can assume $l'=l+2$, and $(a_{l+1},a_{l+2})\in R$. Then $a_{l+1}=y$ and $a_{l+2}=x$, hence $(y,x)\in R\subseteq S$. It follows that $x=y$, as $S$ is antisymmetric. \par
It follows that the pair $(x,y)$ occurs only once in the sequence~\ref{sequence}. Hence $(a,x)\in R$ and $(y,b)\in R$, i.e. $(a,b)\leq_{R^{(2)}}(x,y)$. By symmetry of the roles of $(x,y)$ and $(a,b)$, it follows that $(a,b)=(x,y)$. Thus $S-R$ consists of a single pair $(x,y)$. Moreover if $z\in X$ is such that $(x,z)$ and $(z,y)$ are both in $S$, then either $z=y$ or $(x,z)\neq (x,y)$, thus $(x,z)\in R$. Similarly, either $z=x$ or $(z,y)\in R$. Since $(x,y)\notin R$, one of the pairs $(x,z)$ or $(z,y)$ is not in $R$, hence $x=z$ or $z=y$. Hence $(x,y)$ is adjacent for $S$. \par
Moreover if $(x',y')\in (X\times X)-(R\cup R\op)$ and $(x',y')\leq_{R^{(2)}}(x,y)$, then $(x',x)\in R$ and $(y,y')\in R$. Since $R\subset S$ and $(x,y)\in S$, it follows that $(x',y')\in S$. If $(x',y')\neq (x,y)$, it follows that $(x',y')\in R$. Hence $(x,y)$ is a minimal element of $(X\times X)-(R\cup R\op)$ for the relation $R^{(2)}$, and this completes the proof of Assertion~2.\par
For Assertion 3, let $(x,y)$ be a minimal element of $(X\times X)-(R\cup R\op)$ for the relation $R^{(2)}$. Then as in the proof of Assertion 2, the transitive closure $S=\sur{R\sqcup\{(x,y)\}}$ consists of the union of $R$ with the set of pairs $(a,b)\in X\times X$ such that $(a,x)\in R$ and $(y,b)\in R$, i.e. $(a,b)\leq_{R^{(2)}}(x,y)$. Since $(x,y)$ is minimal for  $R^{(2)}$, it follows that $S=R\sqcup\{(x,y)\}$. This relation $S$ is clearly antisymmetric if $(x,y)\notin R\cup R\op$, hence $S$ is an order and Assertion~3 follows.\findemo
\npar By analogy with the case of subgroups of a group, we set:
\begin{mth}{Definition} Let $S$ be an order on $X$. The {\rm Frattini subrelation} $\Phi(S)$ of $S$ is defined as the intersection of all the maximal order subrelations of $S$.
\end{mth}
\begin{mth}{Proposition} Let $S$ be an order on $X$. Then 
$$\Phi(S)=\Delta\cup (S-\Delta)^2\mpoint$$
In other words, a pair $(x,y)$ of $X\times X$ is in $\Phi(S)$ if and only either $x=y$, or $x<_Sy$ and $]x,y[_S\neq\emptyset$.
\end{mth}
\pf By Lemma~\ref{maximal}, the maximal order subrelations of $S$ are the relations $S-\{(x,y)\}$, where $(x,y)$ is adjacent for $S$. It follows that $\Phi(S)$ consists of the difference $S-\mathcal{M}_S$. The proposition follows.\findemo
\section{Intervals}
\begin{mth}{Theorem} \label{interval and frattini}Let $R\subseteq S$ be two orders on $X$.
\begin{enumerate}
\item If $\Phi(S)\not\subseteq R$, i.e. if there exist $x,y,z\in X$ such that $x<_Sy<_Sz$ and $x\not\leq_Rz$, then the poset $]R,S[_\Omega$ is contractible.
\item If $\Phi(S)\subseteq R$, then any subset of $S$ containing $R$ is an order. In particular, the poset $]R,S[_\Omega$ is isomorphic to the poset of proper non empty subsets of $S-R$, and it has the homotopy type of a sphere of dimension $|S-R-2|$.
\end{enumerate}
\end{mth}
\pf If $T\in ]R,S[_\Omega$, then there is a maximal element $U$ of $]R,S[_\Omega$ which contains $T$, hence $T\cup\Phi(S)$. Thus $\sur{T\cup\Phi(S)}\in ]R,S[_\Omega$. Now if $\Phi(S)\not\subseteq R$, the maps of posets $T\mapsto \sur{T\cup\Phi(S)}\mapsto \sur{R\cup\Phi(S)}$ show that $]R,S[_\Omega$ is {\em conically contractible}, in the sense of \cite{quillen}~1.5.\par
If $\Phi(S)\subseteq R$, then $R=S-A$, where $A$ is a subset of $\mathcal{M}_S$. The map $\delta:T\mapsto T-R$ is a map of posets from $]R,S[_\Omega$ to the set $]\emptyset,A[$ of proper non empty subsets of $A$, ordered by inclusion. \par
{\bf Claim:} for any subset $B\in]\emptyset,A[$, the relation $R\cup B$ is an order on $X$: indeed $R\cup B=S-C$, where $C=A-B$. The case $|C|=1$ is Assertion~1 of Lemma~\ref{maximal}. The general case of the claim follows by induction on $|C|$, from the fact that if $T\subseteq S$, then $\mathcal{M}_S\cap T\subseteq \mathcal{M}_T$.\par
Now the map $\upsilon:B\mapsto R\cup B$ is a map of posets from $]\emptyset,A[$ to $]R,S[_\Omega$. Clearly the maps $\delta$ and $\upsilon$ are inverse isomorphisms of posets, and the poset of proper non empty subsets of a finite set of cardinality $c$ has the homotopy type of the sphere $S^{c-2}$: indeed,  its geometric realization is the boundary of the standard $(c-1)$-simplex in $\R^c$, consisting of the points $(x_1,\ldots,x_c)\in\R^c$ such that $x_i\geq 0$ and $\sum_{i=1}^c\limits x_i=1$.\findemo
\begin{mth}{Corollary} Let $R\subseteq S$ be orders on $X$. Then the M\"obius function $\mu_{\Omega }(R,S)$ of the poset $\Omega $ is equal to 0 if $\Phi(S)\not\subseteq R$, and to $(-1)^{|S-R|}$ otherwise.\findemo
\end{mth}
\pf Indeed, the M\"obius function $\mu_{\Omega }(R,S)$ is equal to the reduced Euler-Poincar\'e characteristic of the poset $]R,S[_\Omega$.\findemo
\section{Upper intervals}
Let $R$ be an order on $X$. This section deals with the homotopy type of the poset $]R,\,.\,[_\Omega$.\par
\begin{mth}{Lemma} \label{lemmamini}Let $R$ be an order on $X$, and $(a,b)\in (X\times X)$. Then $(a,b)$ is a minimal element of $(X\times X)-(R\cup R\op)$ for the relation $R^{(2)}$ if and only if $(a,b)\notin R\cup R\op$ and
\begin{equation}\label{mini}
{]\,.\,,a[_R}\, \subseteq\, {]\,.\,,b[_R} \;\hbox{and}\;{]b,\,.\,[_R}\, \subseteq\, {]a,\,.\,[_R}\mpoint
\end{equation}
\end{mth}
\pf If  $(a,b)$ is a minimal element of $(X\times X)-(R\cup R\op)$ for the relation $R^{(2)}$, then $(a,b)\notin R\cup R\op$. If $x\in X$ and $x<_Ra$, then $(x,b)<_{R^{(2)}}(a,b)$. Hence $(x,b)\in R\cup R\op$. If $(x,b)\in R\op$, then $b\leq_Rx<_Ra$, hence $(b,a)\in R$, a contradiction. Hence $(x,b)\in R$, i.e. $x\leq_Rb$. Moreover $x\neq b$ since $(b,a)\notin R$. Hence $x<_Rb$, i.e. $x\in {]\,.\,,b[_R}$. This shows that ${]\,.\,,a[_R}\, \subseteq\, {]\,.\,,b[_R}$. Similarly, if $y\in X$ and $b<_Ry$, then $(a,y)<_{R^{(2)}}(a,b)$, so $(a,y)\in R\cup R\op$. If $(a,y)\in R\op$, then $b<_Ry\leq_R a$, and $(b,a)\in R$, a contradiction. Hence $(a,y)\in R$, i.e. $a\leq_Ry$. Moreover $a\neq y$ since $b\not\leq_Ra$. Hence $y\in {]a,\,.\,[_R}$, thus ${]b,\,.\,[_R}\, \subseteq\, {]a,\,.\,[_R}\}$, and Condition~\ref{mini} holds.\par
Conversely, if $(a,b)\notin R\cup R\op$, and if Condition~\ref{mini} holds, suppose that $(x,y)\in X\times X$ and $(x,y)\leq_{R^{(2)}}(a,b)$, i.e. $x\leq_Ra$ and $b\leq_Ry$. If $x\neq a$, then $x\in {]\,.\,,a[_R}\;\subseteq\,{]\,.\,,b[_R}$, hence $x<_Rb\leq_Ry$, thus $(x,y)\in R$. Similarly if $b\neq y$, then $y\in {]b,\,.\,[_R}\,\subseteq\,{]a,\,.\,[_R}$, hence $x\leq_Ra<_Ry$, and $(x,y)\in R$. Hence if $(x,y)<_{R^{(2)}}(a,b)$, then $(x,y)\in R\subseteq R\cup R\op$, so $(a,b)$ is a minimal element of $(X\times X)-(R\cup R\op)$ for the relation $R^{(2)}$.\findemo
This motivates the following notation:
\pagebreak[3]
\begin{mth}{Notation}
Let $R\in\Omega $. Set
$$\mathcal{E}_R=\{(a,b)\in (X\times X)-(R\cup R\op)\bmid {]\,.\,,a[_R}\, \subseteq\, {]\,.\,,b[_R} \;\hbox{and}\;{]b,\,.\,[_R}\, \subseteq\, {]a,\,.\,[_R}\}\mpoint$$
\end{mth}
\begin{mth}{Theorem} \label{upper intervals}Let $R$ be an order on $X$.
\begin{enumerate}
\item If there exists $(a,b)\in\mathcal{E}_R$ such that $(b,a)\notin \mathcal{E}_R$, then $]R,\,.\,[_\Omega$ is contractible.
\item Otherwise $\mathcal{E}_R=(\mathcal{E}_R)\op$, and $E_R=\Delta\sqcup\mathcal{E}_R$ is an equivalence relation on~$X$, which can be defined by 
$$(a,b)\in E_R\;\Leftrightarrow\;\left\{\begin{array}{l}{]\,.\,,a[_R}={]\,.\,,b[_R}\\ {]a,\,.\,[_R}={]b,\,.\,[_R}\end{array}\right.\mpoint$$
If $X_1,\ldots,X_r$ are the equivalence classes for this relation, let $\Omega_i=\Omega_{X_i}$ and $\Delta_i=\Delta_{X_i}$, for $i\in\{1,\ldots,r\}$. Then there is a homotopy equivalence
$$]R,\,.\,[_\Omega\,\cong\,]\Delta_1,\,.\,[_{\Omega_1}\,*\,]\Delta_2,\,.\,[_{\Omega_2}\,*\,\ldots\,*\,]\Delta_r,\,.\,[_{\Omega_r}\mvirg$$
where $P*Q$ in the right hand side denotes the {\em join} of two posets $P$ and~$Q$ (cf. \cite{quillen}~Proposition~1.9).
\end{enumerate}
\end{mth}
\pf 
If $(a,b)\in\mathcal{E}_R$, then $(a,b)$ is a minimal element of $(X\times X)-(R\cup R\op)$ for the relation $R^{(2)}$, by Lemma~\ref{lemmamini}. It follows from Lemma~\ref{maximal} that the relation $S=R\sqcup\{(a,b)\}$ is an order on $X$. The poset $[S,\,.\,[_\Omega$ has a smallest element $S$, hence it is a contractible subposet of $]R,\,.\,[_\Omega$. \par
Let $i:[S,\,.\,[_\Omega\hookrightarrow ]R,\,.\,[_\Omega$ denote the inclusion map. For $T\in ]R,\,.\,[_\Omega$, the poset
$$i_{T}=\{U\in [S,\,.\,[_\Omega\,\bmid U\supseteq T\}$$
is equal to the set of order relations containing $S\cup T$. If $\sur{S\cup T}$ is an order, i.e. if it is antisymmetric, then $i_T$ has a smallest element $\sur{S\cup T}$. Otherwise $i_T$ is empty.\par
Let
$$A_S=\{T\in]R,\,.\,[_\Omega\bmid \sur{S\cup T}\;\hbox{is antisymmetric}\}\mpoint$$
In particular $A_S\supseteq [S,\,.\,[_\Omega$, and the inclusion $i:[S,\,.\,[_\Omega\hookrightarrow A_S$ is a homotopy equivalence (\cite{quillen}, Proposition~1.6). Hence $A_S$ is contractible.\par
Now it is easy to check that $\sur{S\cup T}=\sur{\{(a,b)\}\cup T}$ is antisymmetric if and only if $(b,a)\notin\nolinebreak T$.\par
Let 
$$(A_S)_*=\{T\in A_S\bmid\, ]\,.\,,T[_{A_S}\;\hbox{is not contractible}\}\mpoint$$
Then by standard results (e.g. \cite{homol} or~\cite{benson2}, Proposition~6.6.5), if $B$ is any poset with $(A_S)_*\subseteq B\subseteq A_S$, the inclusion $B\hookrightarrow A_S$ is a homotopy equivalence. Now $]\,.\,,T[_{A_S}=]R,T[_\Omega$, for any $T\in A_S$. Let 
$$B=\{T\in ]R,\,.\,[_\Omega\bmid (b,a)\notin T\;\hbox{and}\;\Phi(T)\subseteq R\}\mpoint$$
Then $(A_S)_*\subseteq B\subseteq A_S$ by Theorem~\ref{interval and frattini}, thus $B$ is contractible.\par
Let $T\in ]R,\,.\,[_\Omega$ such that $\Phi(T)\subseteq R$. If $(b,a)\in T$, then $R\sqcup\{(b,a)\}\subseteq T$ is an order, by Theorem~\ref{interval and frattini}. By Lemma~\ref{maximal}, the pair $(b,a)$ is a minimal element of $(X\times X)-(R\cup R\op)$ for the relation $R^{(2)}$. Hence $(b,a)\in\mathcal{E}_R$, by Lemma~\ref{lemmamini}.\par
Hence if there exist $(a,b)\in\mathcal{E}_R$ such that $(b,a)\notin \mathcal{E}_R$, then $\big(]R,\,.\,[_\Omega\big)_*=B$ is contractible, which proves Assertion~1.\par
Otherwise $\mathcal{E}_R\subseteq \mathcal{E}_R\op$, hence $\mathcal{E}_R= \mathcal{E}_R\op$, and
$$\mathcal{E}_R=\{(a,b)\in(X\times X)-(R\cup R\op)\bmid\; ]\,.\,,a[_R=]\,.\,,b[_R\;\hbox{and}\;]a,\,.\,[_R=]b,\,.\,[_R\}\mpoint$$
Then $E_R=\Delta\sqcup\mathcal{E}_R$ is an equivalence relation: it is obviously reflexive and symmetric. To show that it is transitive, it suffices to show that if $(a,b), (b,c)\in \mathcal{E}_R$, then $(a,c)\in E_R$. The assumption implies that $ ]\,.\,,a[_R=]\,.\,,b[_R=]\,.\,,c[_R$ and $]a,\,.\,[_R=]b,\,.\,[_R=]c,\,.\,[_R$. If $(a,c)\in R$, then either $a=c$ or $a\in]\,.\,,c[_R$. In the latter case $a\in ]\,.\,,a[_R$, a contradiction. Thus $a=c$. Now if $(a,c)\in R\op$, then $c\leq_Ra$. Thus again if $a\neq c$, then $c\in  ]\,.\,,a[_R=]\,.\,,c[_R$, a contradiction. In both cases $a=c$. Thus either $a=c$, or $(a,c)\notin R\cup R\op$, and then $(a,c)\in\mathcal{E}_R$. Hence $E_R$ is transitive, so it is an equivalence relation.\par
Let $X_1,\ldots,X_r$ be the equivalence classes of $E_R$. If $T$ is an order containing $R$, for each $i\in\{1,\ldots,r\}$, denote by $T_{|i}$ the restriction of $T$ to $X_i$, i.e. $T_{|i}=T\cap(X_i\times X_i)$. Then $T_{|i}$ is an order on $X_i$, in other words $T_{|i}\in \Omega_i=\Omega_{X_i}$. Let
$$P=[\Delta_1,\,.\,[_{\Omega_1}\times[\Delta_2,\,.\,[_{\Omega_2}\times\ldots\times[\Delta_r,\,.\,[_{\Omega_r}$$
denote the product poset. The map 
$$e:[R,\,.\,[_{\Omega}\to P$$
defined by $e(T)=(T_{|1},\ldots,T_{|r})$ is a map of posets.\par
{\bf Claim:} Let $(T_1,\ldots,T_r)\in P$. Then $S={R\cup T_1\cup\ldots\cup T_r}$ is an order on~$X$.\par
Indeed $S$ is clearly reflexive. It is also antisymmetric: if $(a,b)\in S\cap S\op$, then we can assume that:
\begin{itemize}
\item either $(a,b)\in R\cap R\op$, and then $a=b$, since $R$ is antisymmetric.
\item or $(a,b)\in R$ and $(b,a)\in T_i-R$, for some $i\in\{1,\ldots,r\}$. Then $a\neq b$, thus $(b,a)\in \mathcal{E}_R$, hence ${]\,.\,,a[_R}={]\,.\,,b[_R}$, and $a\in {]\,.\,,b[_R}={]\,.\,,a[_R}$, a contradiction.
\item or $(a,b)\in T_i-R$ (thus $a\neq b$) and $(b,a)\in T_j$, for some $i,j\in \{1,\ldots,r\}$. In particular $(a,b)\in E_R$, so $a$ and $b$ are in the same equivalence class for $E_R$, hence $i=j$. But then $(a,b)\in T_i\cap T_i\op$, thus $a=b$ since $T_i$ is antisymmetric. This is again a contradiction.
\end{itemize}
Hence $S$ is antisymmetric. It is also transitive, for if $(a,b)\in S-\Delta$ and $(b,c)\in S-\Delta$, then:
\begin{itemize}
\item either $(a,b)\in R$ and $(b,c)\in R$, then $(a,c)\in R\subseteq S$, since $R$ is transitive.
\item or $(a,b)\in R-\Delta$ and $(b,c)\in T_i-R$, for some $i\in\{1,\ldots,r\}$. Then $a<_Rb$, thus $a\in{]\,.\,,b[_R}={]\,.\,,c[_R}$, and then $(a,c)\in R\subseteq S$.
\item or $(a,b)\in T_i-R$ and $(b,c)\in R-\Delta$, for some $i\in\{1,\ldots,r\}$. Then $b<_Rc$, i.e. $c\in]b,\,.\,,[_R=]a,\,.\,,[_R$, hence $(a,c)\in R\subseteq S$.
\item or $(a,b)\in T_i-R$ and $(b,c)\in T_j-R$, for some $i,j\in \{1,\ldots,r\}$. Then $(a,b)\in \mathcal{E}_R$ and $(b,c)\in \mathcal{E}_R$, and $a,b,c$ are in the same equivalence class for $E_R$. Hence $i=j$, thus $a<_{T_i}b<_{T_i}c$, thus $a<_{T_i}c$, and $(a,c)\in T_i\subseteq S$.
\end{itemize}
This completes the proof of the above claim.\par
Now let 
$$f:P\to [R,\,.\,[_\Omega$$
be the map defined by $f(T_1,\ldots,T_r)=R\cup T_1\cup\ldots\cup T_r$. It is clearly a map of posets.\par
Moreover for any $T\in[R,\,.\,[_\Omega$
$$fe(T)= R\cup T_{|1}\cup\ldots\cup T_{|r}\subseteq T\mpoint$$
Conversely, for any $(T_1,\ldots,T_r)\in P$
$$ef(T_1,\ldots,T_r)=e(R\cup T_1\cup\ldots\cup T_r)=(T_1',\ldots,T'_r)\mvirg$$
where $T'_i=(R\cup T_1\cup\ldots\cup T_r)\cap(X_i\times X_i)\supseteq T_i$, for any $i\in\{1,\ldots,k\}$. Actually $T'_i=T_i$, since $R\cap(X_i\times X_i)=\Delta_i$~: indeed $R\cap(X_i\times X_i)\subseteq R\cap E_R=\Delta$. Hence $ef$ is the identity map of the poset $P$. This also shows that $e(R)=(\Delta_1,\ldots,\Delta_r)$.\par
Conversely, suppose that $T\in[R,\,.\,[_\Omega$ is such that $e(T)=(\Delta_1,\ldots,\Delta_r)$. Then $T_{|i}=\Delta_i$ for each $i\in\{1,\ldots,r\}$. If $T\neq R$, let $(a,b)\in T-R$, minimal for the relation $R^{(2)}$.  Then $a\neq b$, and $(b,a)\notin R$, since $R\subseteq T$ and $T$ is antisymmetric. Hence $(a,b)$ is also a minimal element of $(X\times X)-(R\cup R\op)$, thus $(a,b)\in\mathcal{E}_R$. It follows that $(a,b)\in T_i$ for some $i\in\{1,\ldots,r\}$. As $a\neq b$, this contradicts the assumption $T_{|i}=\Delta_i$. Hence $T=R$.\par
Finally, let $(T_1,\ldots,T_r)\in P$ such that $f(T_1,\ldots,T_r)=R$. It follows that 
$(T_1,\ldots,T_r)= ef(T_1,\ldots,T_r)=e(R)=(\Delta_1,\ldots,\Delta_r)$.\par
Now the maps $e$ and $f$ restrict to maps of posets
$$\xymatrix{]R,\,.\,[_\Omega\ar[r]<.5ex>^-e&\big([\Delta_1,\,.\,[_{\Omega_1}\times[\Delta_2,\,.\,[_{\Omega_2}\times\ldots\times[\Delta_r,\,.\,[_{\Omega_r}\big)-\{(\Delta_1,\ldots,\Delta_r)\}\ar[l]<.5ex>^-f}$$
such that $fe\subseteq \Id$ and $ef= \Id$. Hence $e$ and $f$ are inverse homotopy equivalences (\cite{quillen}~1.3). This completes the proof of the Theorem, since the poset$$\big([\Delta_1,\,.\,[_{\Omega_1}\times[\Delta_2,\,.\,[_{\Omega_2}\times\ldots\times[\Delta_r,\,.\,[_{\Omega_r}\big)-\{(\Delta_1,\ldots,\Delta_r)\}$$
is homotopic to the join $]\Delta_1,\,.\,[_{\Omega_1}\,*\,]\Delta_2,\,.\,[_{\Omega_2}\,*\,\ldots \,*\,]\Delta_r,\,.\,[_{\Omega_r}$, by~\cite{quillen}~Proposition~1.9.\findemo
\begin{mth}{Theorem} \label{delta} Let $X$ be a set of cardinality $n$. Then the poset $\Omega_X-\{\Delta_X\}$ of non trivial orders on $X$ has the homotopy type of a sphere of dimension ${n-2}$.
\end{mth}
\pf Set $\Omega^\sharp=\Omega-\{\Delta\}$ (where as above $\Omega=\Omega_X$ and $\Delta=\Delta_X$). If $R\in \Omega^\sharp$, then $]\,.\,,R[_{\Omega^\sharp}=]\Delta,R[_\Omega$. If $\Phi(R)\neq\Delta$, then $]\Delta,R[_\Omega$ is contractible by Theorem~\ref{interval and frattini}. It follows that the inclusion
$$B=\{R\in \Omega^\sharp\bmid\Phi(R)=\Delta\}\hookrightarrow \Omega^\sharp$$
is a homotopy equivalence.\par
Now saying that $\Phi(R)=\Delta$ is equivalent to saying that there is no chain $x<_Ry<_Rz$ in $X$. If $R\in B$, let
\begin{eqnarray*}
R_-&=&\{x\in X\mid\exists y\in X,\;x<_Ry\}\\
R_+&=&\{y\in X\mid\exists x\in X,\;x<_Ry\}\\
R_0&=&X-(R_-\cup R_+)\mpoint
\end{eqnarray*}
Define $h(R)=R\cup (R_-\times R_+)$. Then $h(R)$ is an order on $X$: it is obviously reflexive, since if $(a,b)\in h(R)\cap h(R)\op$, then:
\begin{itemize}
\item either $(a,b)\in R$ and $(b,a)\in R$, hence $a=b$.
\item or $(a,b)\in R$ and $(b,a)\in R_-\times R_+$. Then $b\in R_-$ and $a\in R_+$. Moreover $a\neq b$, and $a<_Rb$. Thus $a\in R_-$ and $b\in R_+$. This is a contradiction since $R_-\cap R_+=\emptyset$.
\item or $(b,a)\in R$ and $(a,b)\in R_-\times R_+$. The argument of the previous case applies to the pair $(b,a)$, and this gives a contradiction again.
\item or $(a,b)\in R_-\times R_+$ and $(b,a)\in R_-\times R_+$. Again, since $R_-\cap R_+=\emptyset$, this is a contradiction.
\end{itemize}
Finally if $(a,b)\in h(R)$ and $(b,c)\in h(R)$, then either $a=b$ or $b=c$, so $h(R)$ is trivially transitive. It follows that $h(R)\in\Omega^\sharp$. Clearly $R\subseteq h(R)$.\par
If $R\subseteq S\in B$, then clearly $R_-\subseteq S_-$ and $R_+\subseteq S_+$, thus $h(R)\subseteq h(S)$. Moreover $h(R)_-=R_-$ and $h(R)_+=R_+$, thus $h^2=h$. It follows that $h$ induces a homotopy equivalence 
$$B\to h(B)=\{R\in B\mid h(R)=R\}\mvirg$$
inverse to the inclusion $h(B)\hookrightarrow B$. Here $h(B)$ is viewed as a subposet of $B$, i.e. $h(B)$ is ordered by inclusion.\par
Now if $R\in B$ and $h(R)=R$, then the pair $\pi(R)=(R_-,R_+)$ is a pair $(V,U)$ of disjoint non empty subsets of $X$. Conversely, such a pair $(V,U)$ yields an element $\omega(V,U)$ of $h(B)$ defined by $\omega(V,U)=\Delta\cup (V\times U)$. The maps $\pi$ and $\omega$ induce inverse isomorphisms of posets between $h(B)$ and the poset
$$D_2(X)=\{(V,U)\mid \emptyset\neq V\subset X,\;\emptyset\neq U\subset X,\;V\cap U=\emptyset\}\mvirg$$
ordered by the product order $(V',U')\subseteq (V,U) \Leftrightarrow V'\subseteq V$ and $U'\subseteq U$ (note that if $U$ and $V$ are disjoint non empty subsets of $X$, then $U$ and $V$ are {\em proper} subsets of $X$). \par
Let $D_1(X)$ denote the set of non empty proper subsets of $X$, ordered by inclusion. The projection map $g:D_2(X)\to D_1(X)$ defined by $g(V,U)=U$ is a map of posets. Moreover for $W\in D_1(X)$
$$g_W=\{(V,U)\in D_2(X)\mid g(V,U)\supseteq W\}=\{(V,U)\in D_2(X)\mid U\supseteq W\}$$
is contractible via the maps $(V,U)\mapsto (V,W)\subseteq (X-W,W)$. It follows that $g$ is a homotopy equivalence (\cite{quillen}, Proposition 1.6). Hence the posets $\Omega^\sharp$ and $D_1(X)$ are homotopy equivalent. \par
This completes the proof, since if $X$ has cardinality~$n$, then $D_1(X)$ has the homotopy type of a sphere of dimension $n-2$.\findemo
\begin{mth}{Corollary} Let $R$ be an order on $X$. Then:
\begin{enumerate}
\item if $\mathcal{E}_R\neq \mathcal{E}_R\op$, the poset $]R,\,.\,[_\Omega$ is contractible.
\item if $\mathcal{E}_R= \mathcal{E}_R\op$, then $E_R=\Delta\sqcup\mathcal{E}_R$ is an equivalence relation on $X$. The poset $]R,\,.\,[_\Omega$ is homotopy equivalent to a sphere of dimension $n-r-1$, where $r$ is the number of equivalence classes of $E_R$.
\end{enumerate}
\end{mth}
\pf Assertion 1 follows from Theorem~\ref{upper intervals}. For Assertion~2, by Theorem~\ref{upper intervals}, there is a homotopy equivalence
$$]R,\,.\,[_\Omega\,\cong\,]\Delta_1,\,.\,[_{\Omega_1}\,*\,]\Delta_2,\,.\,[_{\Omega_2}\,*\,\ldots\,*\,]\Delta_r,\,.\,[_{\Omega_r}\mvirg$$
where $\Omega_i=\Omega_{X_i}$ is the poset of posets on the equivalence class $X_i$ of $X$ for the relation $E_R$.\par
By Theorem~\ref{delta}, each poset $]\Delta_i,\,.\,[_{\Omega_i}$ has the homotopy type of a sphere of dimension $|X_i|-2$. Since for $d,e\in\N$, the join $S^d*S^e$ is homotopy equivalent to $S^{d+e+1}$ (where $S^d$ denotes a sphere of dimension $d$), it follows that $]R,\,.\,[_\Omega$ has the homotopy type of a sphere of dimension
$$(|X_1|-2)+(|X_2|-2)+\ldots+(|X_r|-2)+r-1=n-2r+r-1=n-r-1\mvirg$$
as was to be shown.\findemo

\centerline{\rule{5ex}{.1ex}}
\begin{flushleft}
Serge Bouc\\
CNRS-LAMFA, Universit\'e de Picardie\\
33 rue St Leu, 80039 - Amiens - France\\
{\tt email: serge.bouc@u-picardie.fr}\\
{\tt web: http://www.lamfa.u-picardie.fr/{$\sim$}bouc}
\end{flushleft}
\end{document}